# A multi-criterion simulation model to determine dengue outbreaks

Piotr Jakubowski, Hasitha Erandi, Anuradha Mahasinghe, Sanjeewa Perera and Andrzej Ameljańczyk

*Abstract*: *In this study, we develop a multi-criteria model to identify dengue outbreak periods. To validate the model, we performed a simulation using dengue transmission-related data in Sri Lanka's Western Province. Our results indicated that the developed model can be used to predict a dengue outbreak situation in a given region up to one month in advance.*

## 1. INTRODUCTION

Dengue is a vector-borne viral infection that has recently become a critical public health problem [6, 15]. Geographically, the disease has been closely associated with tropical and subtropical climatic regions [15]. Global warming and spatial changes in climate factors may expand the suitable areas for vector habitats and affect the biology and ecology of the vectors [15], increasing the risk of disease transmission. In addition to climate change, changes in other global drivers such as urbanization and human mobility have also increased the rate of disease transmission.

Recently, a significant rise in the number of dengue cases was reported in Sri Lanka. According to the data released by Sri Lanka's Ministry of Health, the worst outbreak was reported in 2017, and the highest number of dengue cases was reported in July 2017 (with 46.5% from the Western Province). Figure 1 illustrates the reported dengue cases of the Western Province of the country from January 2010 to April 2020.

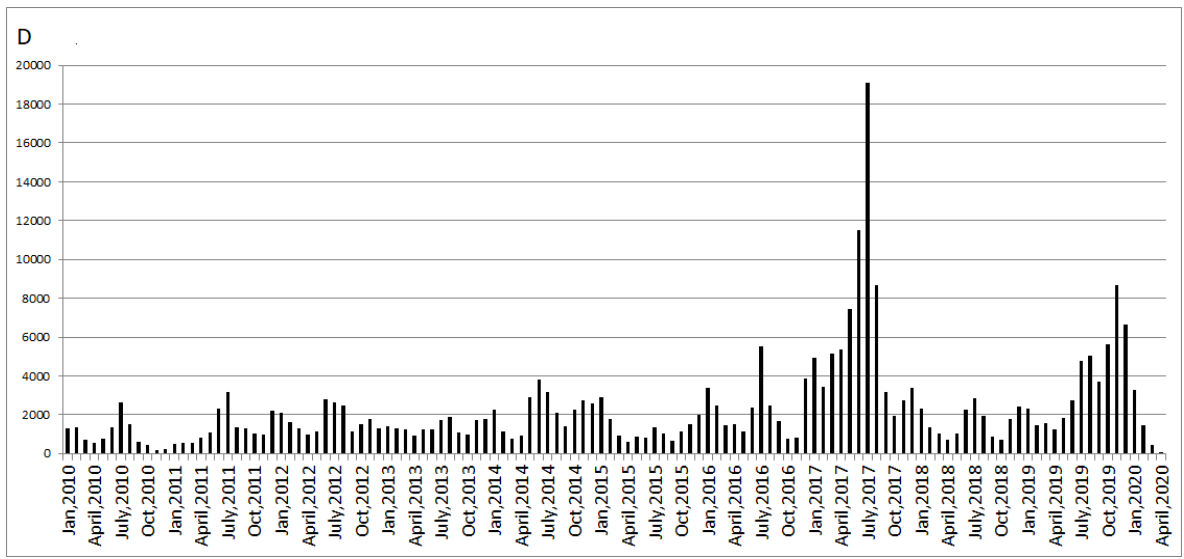

**Figure 1.** Reported dengue cases from January 2010 to April 2020, Western Province, Sri Lanka (Source: Epidemiology Unit, Ministry of Health, Sri Lanka)

Although the first dengue vaccine was licensed in 2015, the vaccine's performance is dependent on its serostatus; the development of the vaccine is still at the stage of experimentation [28]. Hence, the main prevention strategy is vector controlling. For countries with limited resources (like Sri Lanka), the limitation of resources becomes a major obstacle to the control process. Therefore, predicting future outbreaks provides opportunities to plan efficient control strategies.

Several dengue outbreak prediction models have been proposed by researchers and are applicable in different contexts [9, 11, 17, 18, 25]. Although a number of works have considered climate change and/or human travel patterns when identifying dengue outbreaks [11,25], the applicability of these models is confined to those countries with relevant seasonality and mobility patterns. This implies that there is still a need for new models that take the factors that affect dengue transmission in local regions into account.

In order to capture the regional potential risk of dengue transmission in Sri Lanka, we consider regional climate analysis data and inter-regional mobility pattern. On the other hand, regional variation of the risk of dengue transmission depends on population density. To bring this factor into model, the susceptible and infected population densities are considered.

In order to formulate the relevant risk functions, we use techniques based on fuzzy set theory [26, 27]. The membership functions for monthly rainfall, average temperature, relative humidity, and inter-regional mobility pattern are constructed in order to capture the potential risk from each factor. These membership functions are constructed by utilizing existing models in the literature [7, 14, 21, 24]. Subsequently, we calculate the combined effect of the membership functions using the correlation coefficient. Then, we formulate a function to measure the closeness similarity of both the potential and variation risk functions to the ideal situation of dengue transmission. Consequently, the problem turns out to be a two-dimensional optimization problem for which a Pareto optimum provides a realistic solution.

We apply the solution technique to data sets on dengue incidence, climate factors, and human mobility in the Western Province of Sri Lanka from 2010 to 2018 (first to validate the model, and then to predict the dengue outbreak periods). Then, we compare our results with what could be generated by applying a generalized linear regression on the actual data.

The remainder of the paper is organized as follows: First, we discuss the study area in detail in Section 2. Then, we define regional potential risk and local regional variation of dengue transmission in Section 3. In the same section, we describe the Pareto optimization process as well as the algorithm. This is followed by the computational results that we obtained for dengue transmission in the Western Province of Sri Lanka that will be described in Section 4. Our concluding remarks can be found in Section 5.

## 2. STUDY AREA AND DATA

The Western Province is the most densely populated area in the wet zone of the country, where 5,821,710 people live in an area of 3,593 km2 (according to a 2012 census). The province consists of three districts; Colombo, Gampaha, and Kalutara. Furthermore, the province is home to Sri Jayewardenepura (the country's legislative capital) and Colombo (its administrative and business center). Table 1 represents the population, area size, population density, and reported dengue data as a percentage of total dengue cases in Sri Lanka, the Western Province, and its relevant districts in 2018.

Moreover, almost all of the premier educational institutions and the largest number of schools in the country are located in the Western Province, and the public transport service connects the province to the other major cities on the island. Hence, viral disease can be easily transmitted from the Western Province to other provinces due to the high human mobility rates. Figure 2 demonstrates the boundaries of the province.

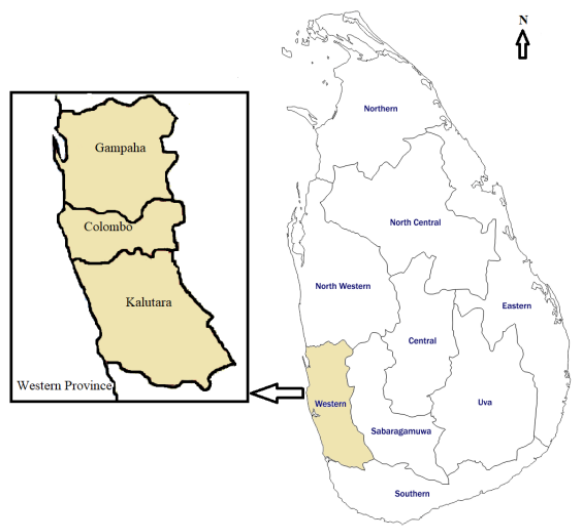

**Figure 2.** Area map of Western province, Sri Lanka

**Table 1:** Population, area, population density and reported dengue data as a percentage of total dengue cases

| Region | Population (million) | Area (km²) | Population. density | Reported dengue data | |
|---|---|---|---|---|---|
| | | | | *2017* | *2018* |
| Sri Lanka | 21M | 65,610 | 327 | 186101 | 51659 |
| Western province | 5.8M | 3,593 | 1,628 | 41% | 37% |
| Colombo district | 2.3M | 676 | 3,438 | 18% | 20% |
| Gampaha district | 2.3M | 1,341 | 1,719 | 17% | 11% |
| Kalutara district | 1.2M | 1,576 | 775 | 5.9% | 6% |

The Western Province is influenced by two monsoon seasons: the Northeast monsoon season (from December through February), and the Southwest monsoon season (from May through September). The weather conditions of the province are generally warm and humid, with a maximum temperature of 33℃ during the day and a. minimum temperature of 22℃ at night. During the period of 2010 through 2015, the average monthly rainfall was 3.2–794.8 mm, and the air humidity fell within a range of 62–95%. Therefore, the Western Province is an ideal region for our study.

For this study, we considered the monthly reported dengue data gain from the Epidemiology Unit, Department of Health, Sri Lanka as well as the monthly rainfall, average temperature, and maximum humidity data for the Western Province from January 2010 through December 2018 from the Department of Meteorology, Sri Lanka. To measure the inter-provincial human mobility, we considered the inter-provincial bus route data from the National Transport Commission, Sri Lanka.

### 3. MODEL DEVELOPMENT

In order to develop the multi-criteria simulation model, we follow the following simple approach:

• Step 1: Define membership functions for each climate factor and mobility factor.

• Step 2: Define regional potential risk function of dengue using combined effect of membership functions.

• Step 3: Define local regional variation using regional susceptible and infected population.

• Step 4: Formulate function to identify closeness similarity to ideal situation of dengue transmission.

• Step 5: Define algorithm and validate model.

*Fuzzy set theory*

The concept of fuzzy set theory was introduced by Lotfi Zadeh in 1965 and has been applied to many fields, including control theory [30]. Consider a non-empty set $U$ and an element $x \in U$. A fuzzy set is defined as a non-empty subset of $U$ U, and function $F: U \to [0,1]$ is called the membership function of the fuzzy set. Fuzzification is the process of assigning the numerical input of a system to fuzzy sets with some degree of membership within an interval of [0,1]. If the membership value is 0, then $x$ does not belong to a given fuzzy set; if the membership value is 1, then $x$ completely belongs within the fuzzy set. Any value between 0 and 1 represents the degree of uncertainty that the value belongs in the set [10]. Any value between 0 and 1 represents the degree of uncertainty that the value belongs to the set [19]. Note that defining the potential risk functions of dengue transmission is an uncertain process with insufficiently reliable data sources. Moreover, the influence of climate factors on dengue transmission is inter-relational [16, 21]. Therefore, we adopt fuzzy set theory and its tools to capture the combined effect on the potential risk of dengue transmission.

*Regional potential risk of dengue transmission*

The transmission potential of dengue within a region is influenced by local climate factors such as rainfall data, temperature, and humidity as well as global factors such as inter-regional human mobility [16, 21, 29]. Each climate factor has a potential impact on dengue transmission; to evaluate the potential impact of the temperature and humidity factors on dengue transmission, we consider the laboratory confirmed data in [21].

For an example, laboratory tests have confirmed that the ideal survival temperature for all phases of vector development ranges from 20 to 30°C [21]. Immature vector development rates decrease after 34°C and egg, larvae, and pupae development rates decrease at temperatures of lower than 8.3°C. Furthermore, mosquitoes become inactive at temperatures that are lower than 15°C or higher than 36°C. Moreover, the mortality rates of adult mosquitoes increase with increasing temperature above 30°C [14]. Let $X_T(i,t)$ denote the average temperature data for the $i$ th region at time $t$ and $Y_T(i,t)$ (the membership function of dengue transmission subject to average temperature). Considering the above-mentioned laboratory-confirmed data, a function for $Y_T(i,t)$ is defined as in Figure 3, where 0 represents the lowest potential risk and 1 is associated with the highest potential risk.

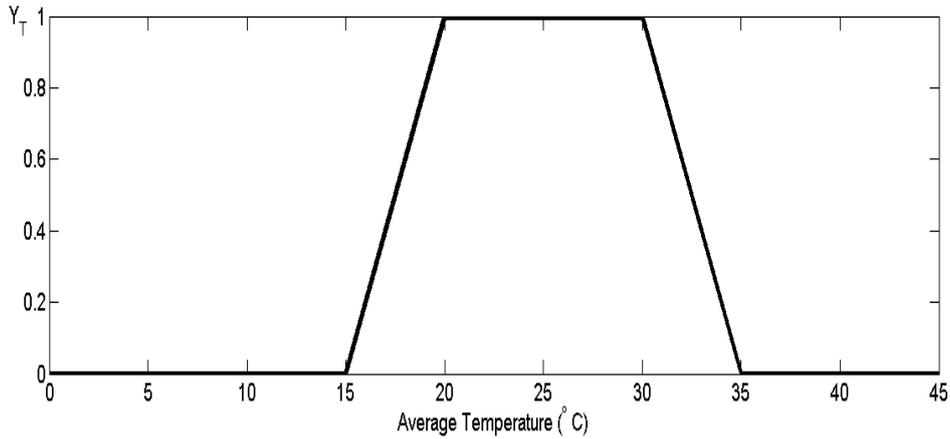

**Figure 3.** Membership function for average temperature

Similarly, membership function $(Y_H(i,t))$ of dengue transmission for humidity data $(X_H(i,t))$ is defined based on laboratory confirmed data in [7] (see Figure 4).

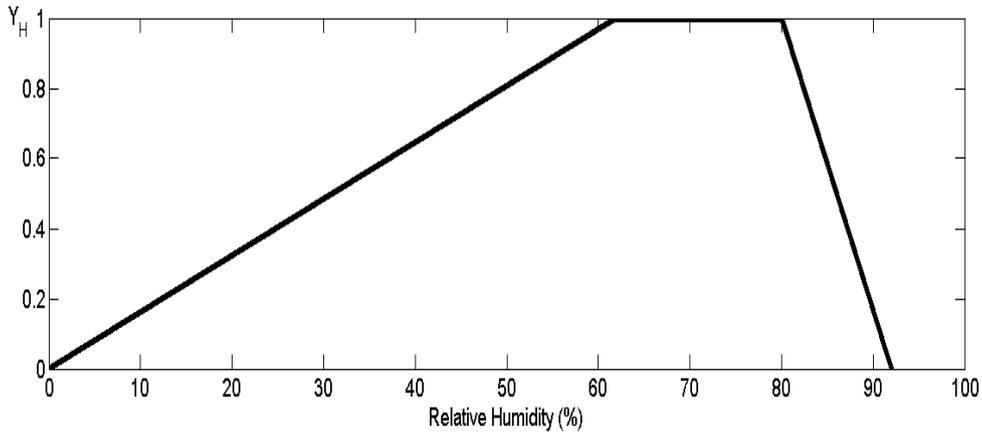

**Figure 4.** Membership function for relative humidity

Rainfall mainly affects dengue by generating vector breeding sites [24]. Since rainfall intensity has both positive and negative effects on dengue breeding sites [24], to evaluate the potential impact of rainfall on dengue transmission, we calculate the minimum and maximum cutoff points for rainfall data $(X_R(i,t))$ that gives the highest correlation value with dengue using Pearson correlation coefficient. Then, the membership function to identify the potential risk functions of dengue transmission for rainfall data $(Y_R(i,t))$ is defined based on the correlation values between the dengue data and rainfall data (see Figure 5).

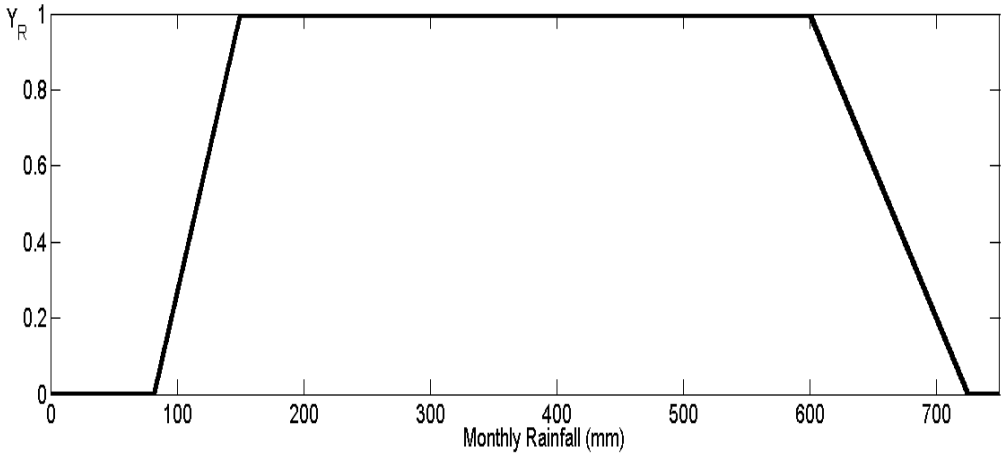

**Figure 5.** Membership function for monthly rainfall

Sri Lanka is a small country with ample human mobility. Therefore, the regional potential risk of dengue transmission depends on its inter-regional human mobility (see Figure 6). The term $X_m(i,j)$ is a mobility–based interaction measure between $i$th and $j$th regions. $I(j,t)$ is infected population density in $j$th region at time $t$. Then, the potential risk of dengue transmission upon mobility can be defined as follows;

$$R_{mob}(i,t) = \sum_{i \sim j} X_m(i,j) I(j,t) . \qquad (1)$$

Notice that, transmission potential of dengue is proportional to the inter-regional mobility and the infected population density within the connected regions. Hence, we define a linear function that range between 0 and 1 with a positive slope as a function for the mobility factor and is denoted by $Y_m(i,t)$, where $C$ denotes the maximum value of $R_{mob}(i,t)$.

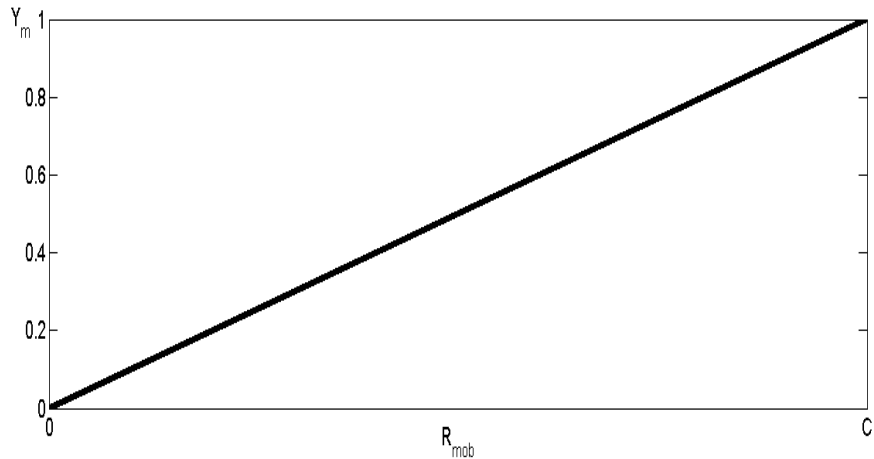

**Figure 6.** Potential risk function for mobility factor

Combining all the membership functions, regional potential risk of dengue transmission $R(i,t)$ can be defined as follows:

$$R(i,t) = Y_r^{c_1} \times Y_T^{c_2} \times Y_H^{c_3} \times Y_m^{c_4}. \qquad (2)$$

Here, $c_1$, $c_2$, $c_3$ and $c_4$ are power of each potential risk function; the influence of each function to the regional potential risk of dengue transmission can be changed by modifying $c_1$, $c_2$, $c_3$ and $c_4$. To identify $c_1$, $c_2$, $c_3$ and $c_4$, a cross-correlation coefficient between climate factor and dengue data has been used.

Notice that, *Aedes* mosquitoes proceed life cycle from eggs to adult through larvae and pupae; their life cycle takes approximately 1–2 weeks or longer depending on the temperature as well as the availability of water and nutrients [5]. Furthermore, the incubation period for an infected human ranges from 3 to 14 days; an infected human experiences a shorter incubation period for dengue viruses with high temperature and favorable humidity [5].

Therefore, it is possible to obtain Equation (3) by including the time lag into rainfall data, temperature data, humidity data and mobility data in Equation (2).

$$R(i,t) = Y_r(i,(t-lr))^{c_1} \times Y_T(i,(t-lT))^{c_2} \times Y_H(i,(t-lH))^{c_3} \times Y_m i,(t-lm)^{c_4}. \qquad (3)$$

Since all the function values for $Y_r, Y_T, Y_H$ and $Y_m$ range within a scale of 0 to 1, the functional values of $R(i,t)$ are in between 0 to 1. Therefore, the 0 value of $R(i,t)$ is associated with the lowest regional potential risk value, and 1 is associated with the highest regional potential risk value.

*Local regional variation of dengue*

Although the environment within a region might be favorable for disease transmission, the speed of disease spread is governed by the susceptible and infected hosts within that region. For example, if all of the population in the $i$th region is infected with dengue at time $t$, then there is no regional variation of the disease at time $t+1$. Therefore, we define a function for the local regional variation of disease $L(i,t)$ within region $i$ at time $t$ based on infected host density $I$ and susceptible host density $S$ in the particular region at time $t-1$.

$$L(i,t) = Y_S(i,(t-1)) \times Y_I(i,(t-1)), \qquad (4)$$

where, $Y_S(i,(t-1))$ and $Y_I(i,(t-1))$ denote the functions of local regional variation upon $S(i,(t-1))$ and $I(i,(t-1))$, respectively. Moreover, the function values for $Y_S(i,(t-1))$ and $Y_I(i,(t-1))$ range within a scale of 0 to 1. Hence, the functional values of $L(i,t)$ are between 0 and 1, with 1 being associated with the highest risk of local regional variation of dengue transmission. Now, we will investigate how the regional potential risk and local regional variation of dengue transmission could be used to forecast outbreak periods over time.

*Similarity Identification*

To identify the outbreak periods over time for a given region $i^*$, the closeness similarity to the effective regional potential and local regional variation should be calculated. Let $R_{ideal}$ and $L_{ideal}$ denote the most ideal values for the regional potential risk and local regional variation of dengue transmission for given region $i^*$. Then, the closeness in similarity to the ideal risk of dengue transmission for each time period could be defined as follows:

$$d_1(i^*,t) = 1 - \frac{R(i^*,t)}{R_{ideal}}, \qquad (5a)$$

$$d_2(i^*,t) = 1 - \frac{L(i^*,t)}{L_{ideal}}. \qquad (5b)$$

If the considered region at time $t$ has an ideal environment for dengue transmission, then the value of $d_1(i^*,t)$ is close to 0. On the other hand, if the risk of local regional variation has a favorable condition for dengue transmission at time $t$, then the value of $d_2(i^*,t)$ is also close to 0. To identify dengue outbreak periods over time, two objective functions must be satisfied. First, the regional environment should be ideal for dengue transmission. Second, the local regional variation should be ideal for disease transmission. That is, we must identify the non-dominating time points with both function values $d_1(i^*,t)$ and $d_2(i^*,t)$ close to 0. To solve this two-dimensional problem, we use an algorithmic version of the notion of Pareto optimization [4, 13].

*Pareto Optimization Algorithm*

The concept of Pareto efficiency was first introduced by Vilfredo Pareto [10, 23] to describe a phenomenon in economics. Later works investigated the possibility of expanding the technique to engineering design and multi-objective optimization [4,8,13]. The solution set of Pareto optimization is known as the Pareto optimal frontier [8], and a solution in a Pareto optimal set cannot be deemed superior to the others in the set without including preference information to rank the competing attributes.

Figure 7 provides a visualization of the general idea of the Pareto optimal frontier. According to Figure 7, the solution space can be divided into two cluster groups; $y^*(x)$ and $y^*(\bar{x})$. Notice that the points in group $y^*(\bar{x})$ are closer to the ideal point $(0,0)$; hence, it contains solution points in Pareto optimal frontier.

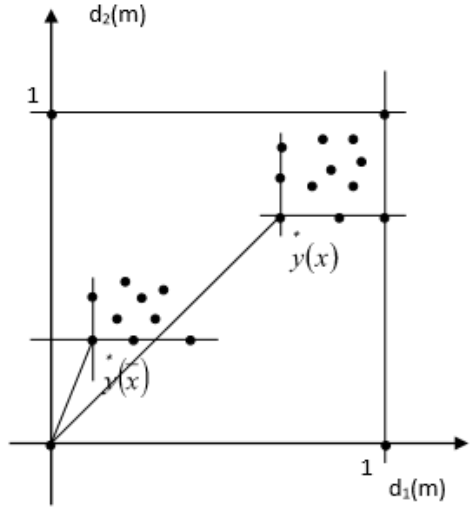

**Figure 7.** Solution space

Recall that our problem has two objective functions (as defined in Equation 5). Hence, for a given region $i^*$, the multi-objective function to be minimized is expressible as follows:

$$Z(i^*, t) = (d_1(i^*, t), d_2(i^*, t)). \qquad (6)$$

It is straightforward to see that first objective function in Equation 6 depends on the regional potential risk of dengue transmission in Equation 3, and the second objective function depends on the local regional variation of dengue in Equation 4. Hence, the constraints related to both objective functions can be expressed as:

$$R(i^*, t) = Y_r(i^*, (t - lr))^{c_1} \times Y_T(i^*, (t - lT))^{c_2} \times Y_H(i^*, (t - lH))^{c_3} \times Y_m(i^*, (t - lm))^{c_4} \qquad (7a)$$

$$L(i^*, t) = Y_S(i^*, (t - 1)) \times Y_I(i^*, (t - 1)) \qquad (7b)$$

Motivated by a previous works [1-3], we developed a multi-objective algorithm to generate computational results for our problem by combining the Pareto optimization technique with the particular objectives. In order to identify the Pareto optimal frontier, we defined a ranking system based on Equation 5 and degree of dominance [1-3]. The algorithm is defined by the following steps.

- Step 1: Develop functions to identify risk of dengue transmission.
- Step 2: Input climate data, mobility data, and population data for each time point.
- Step 3: Calculate rank of each time point based on potential risk and local regional variation and degree of dominance.
- Step 4: Run detection algorithm based on rank values and detection information (refer to Figure 9).
- Step 5: Identify Pareto optimal frontier.

To have a clear understanding about the detection process, a schematic representation of the detection process and the detection algorithm are shown in Figures 8 and 9 respectively. The initial identification in the detection process is done with the help of laboratory-confirmed data. Moreover, the detection-inference mechanism depends on the power of each potential risk function and the local variation risk values. The analysis of the power of each risk function allows us to determine the outbreak time periods that are nearly similar to the actual outbreak time periods.

Notice that the detection process of the model and the algorithm allowed us to classify the obtained results based on real data, confirming the reliability of the results with the model of the occurrence of the outbreak. The reliability of the identified point is a derivative of the distance between the ideal point and the prediction time period. If the distance is small, then the outbreak prediction process is more reliable.

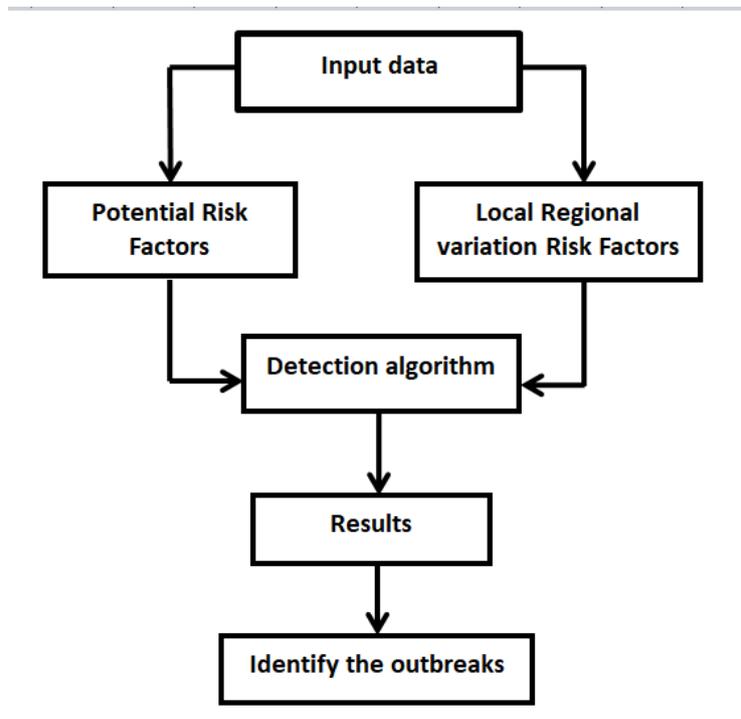

**Figure 8.** Scheme of detection process

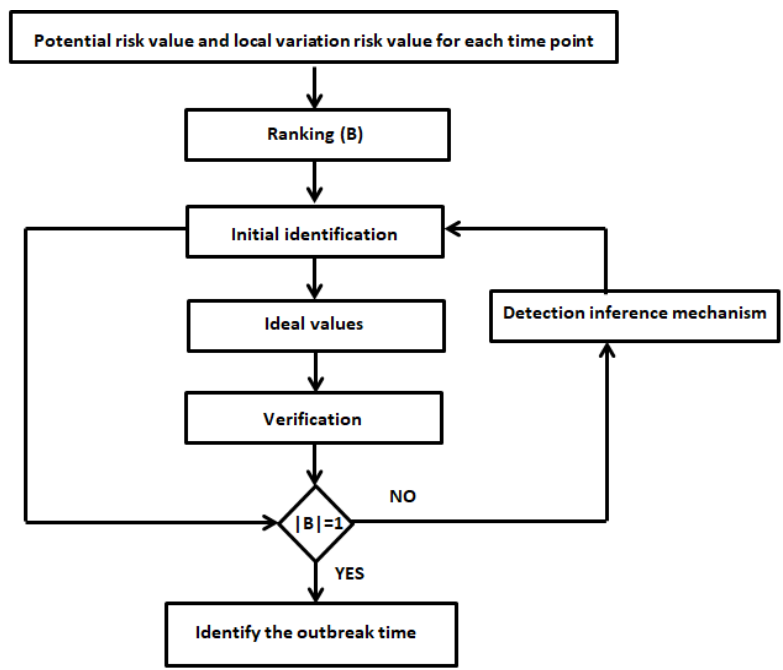

**Figure 9.** Scheme of detection (Pareto simulation) algorithm

## 4. RESULTS AND DISCUSSION

To determine the lag times in Equation 3, we computed the cross-correlation between the dengue data and climate data using Pearson cross-correlation formula. According to the data analysis, the lag times for the highest cross-correlation value between the dengue data and the rainfall, average temperature and maximum humidity are two, three and two months, respectively. For the mobility data in Equation 3, we used one-month-back infected population data.

In order to simulate the model, we used an application simulator. The application was written based on Microsoft .NET Framework technologies. The C# language was used in the implementation process. The design environment that was used to implement the software was Microsoft Visual Studio 2017.

Figure 10 represents the Pareto optimal frontier for the Western province's monthly data from April 2010 through

December, 2018.

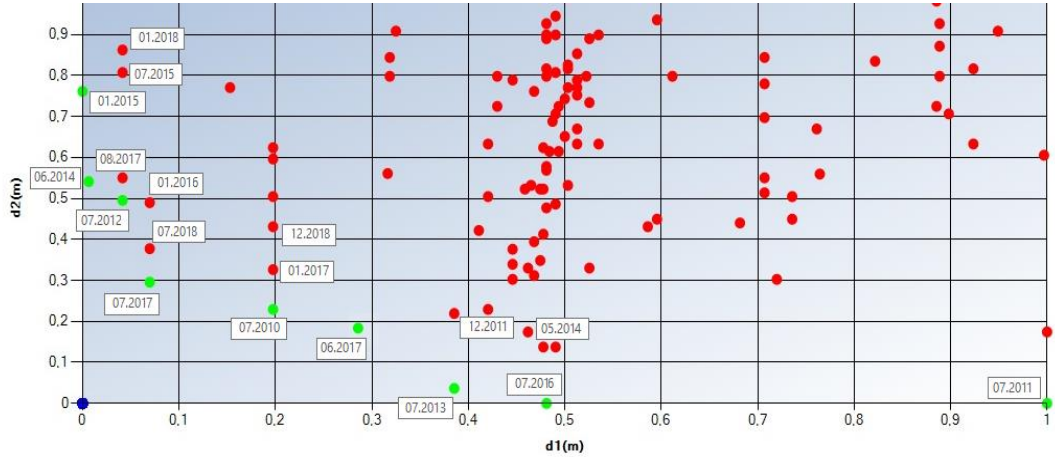

**Figure 10.** Pareto optimal front for Western Province's monthly data from April 2010 through December 2018

According to Figure 10, each July of 2010, 2011, 2012, and 2013, June 2014, January 2015, July 2016, and June and July 2017 are in the Pareto optimal frontier. Comparing the months in the Pareto optimal frontier with the actual dengue data, it can be observed that all of them are months with high dengue incidences. However, according to Figure 2, there were other outbreaks during the considered time period, such as December 2011 and July 2018 (these are not included in the Pareto optimal frontier). Therefore, we considered the close points to the Pareto optimal frontier with similar rank values as illustrated in Figure 11.

Since the proposed model is new to dengue outbreak prediction, we compared the results with generalized linear-regression model [12, 22] based on the same input data. In the generalized linear-regression model, the dependent variable (the number of dengue cases) is assumed to follow the distribution of the dependent parameter only through the linear combination [20].

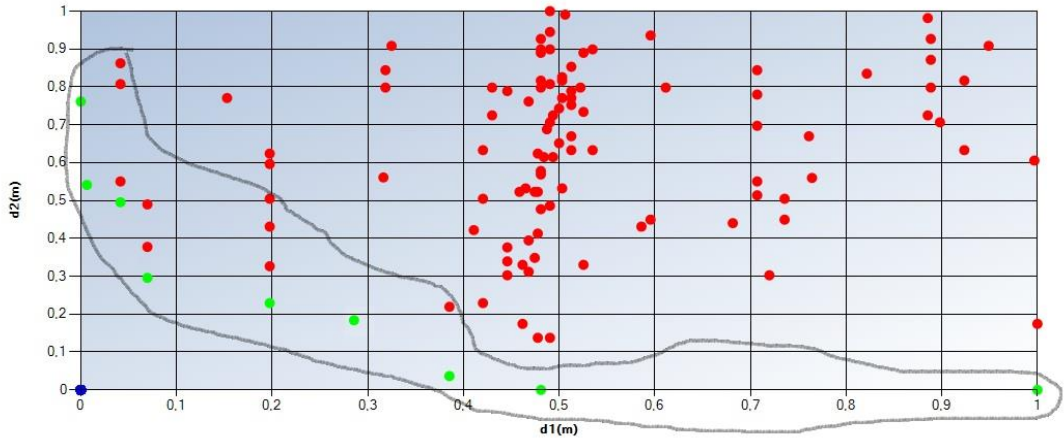

**Figure 11.** Close points to Pareto optimal front for Western Province's monthly data from April 2010 through December 2018

In this study, we used the glmfit package in MATLAB to estimate parameters $b_0$ through $b_6$ based on the distribution fitting with the assumption of linearity $D(t)$, which is expressed as follows:

$$D(i^*, t) = b_0 + b_1. x_R(i^*, (t - lag_R)) + b_2. x_T(i^*, (t - lag_T)) + b_3. x_H(i^*, (t - lag_H)) + b_4. R_{mob}(i^*, (t - lag_{R_m})) + b_5. I(i^*, (t - 1)) + b_6. S(i^*, (t - 1)) \quad (8)$$

Figure 12 represents a comparison of the simulation results and actual data for the linear-regression model.

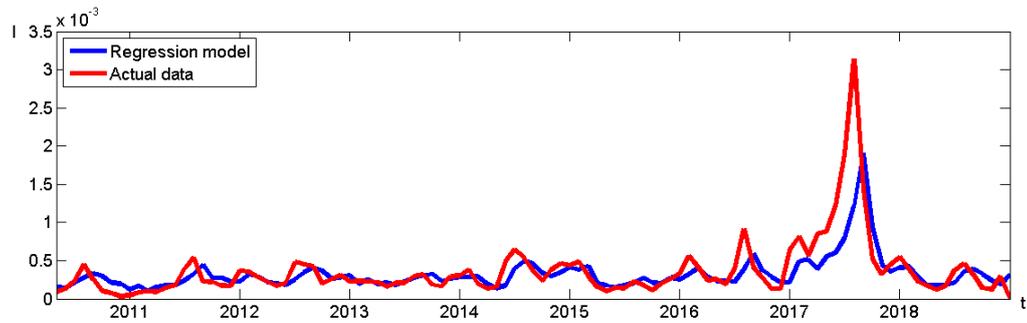

**Figure 12.** Simulation results for generalized linear-regression model for Western Province's monthly data from April 2010 through December 2018 (I is infected population density in Western Province)

Then, we compared the results obtained by the multi-criteria simulation and generalized liner-regression models with the actual dengue outbreaks. Table 2 represents the simulated results for both models and the actual outbreak periods.

**Table 2:** Comparison of results obtained by multi-criteria simulation and generalized liner-regression models and actual dengue outbreaks)

| Actual outbreak | Multi-criterion simulation | Linear-regression |
|---|---|---|
| July 2010 | July 2010 | August 2010 |
| July 2011 | July 2011 | August 2011 |
| December 2011 | December 2011 | January 2012 |
| June 2012 | - | - |
| July 2012 | July 2012 | July 2012 |
| August 2012 | - | - |
| November 2012 | - | November 2012 |
| July 2013 | July 2013 | July 2013 |
| August 2013 | July 2013 | August 2013 |
| January 2014 | - | December 2013 |
| June 2014 | June 2014 | July 2014 |
| November 2014 | - | November 2014 |
| January 2015 | January 2015 | January 2015 |
| January 2016 | January 2016 | February 2016 |
| July 2016 | July 2016 | August 2016 |
| January 2017 | January 2017 | January 2017 |
| May 2017 | May 2017 | - |
| June 2017 | June 2017 | June 2017 |
| July 2017 | July 2017 | July 2017 |
| August 2017 | August 2017 | August 2017 |
| December 2017 | January 2018 | December 2017 |
| July 2018 | July 2018 | August 2018 |
| November 2018 | November 2018 | - |

From Table 2, it can be observed that the Pareto optimal front and the close points contain seven false-positive and false-negative data points. On the other hand, the results of the linear-regression model contain 12 false-positive and false-negative data points. Hence, the false-positive and false-negative answer rate for the multicriteria simulation and linear-regression models are 6.67% and 11.43%, respectively. Therefore, the proposed multi-criteria simulation model is more accurate than the linear-regression model.

Recall that we used rainfall for Equation (3), average temperature, maximum humidity, and mobility data with 2, 3, 2, and a one-month time lag, respectively. The proposed model can be used to identify the outbreak situation of a given region for up to one month using real-time data. Hence, the model can be used as a real-time warning system, helping to provide sufficient time to control the spread of the disease.

**CONCLUSION**

Dengue is a rapidly spreading disease in the world that has a complex transmission mechanism. Since dengue transmission is highly affected by regional climate changes and mobility, estimating future outbreaks is vital for controlling the spread of the disease. Since existing classical models based on parameters are not very sensitive to time nor space, they do not have the full capacity to capture the future outbreak of a given region.

In this work, we have developed a new model to identify a dengue outbreak period in a given region up to one month in advance using real-time rainfall, temperature, humidity mobility, and regional infected and susceptible data.

Though our main concern was identifying the dengue outbreak periods in a given region, the model is also applicable in identifying the source of an epidemic network in given time period; an interesting future task might be to extend the model to identify regions with dengue outbreaks during a given time period.

ACKNOWLEDGEMENT

This work was supported by Grant RPHS/2016/D/05 of the National Science Foundation, Sri Lanka and the European Social Fund under the Operational Program Knowledge Education Development 2014-2020.

**REFERENCES**

[1]. Ameljanczyk A.: Properties of the Algorithm for Determining an Initial Medical Diagnosis Based on a Two-Criteria Similarity Model, *Biuletyn Instytutu Systemów Informatycznych*, vol. 8, pp. 9–16, 2011.
[2]. Ameljanczyk A.: Pareto filter in the process of multi-label classifier synthesis in medical diagnostics support algorithms, *Computer Science and Mathematical Modelling*, vol. 1, pp. 5–10, 2015.
[3]. Ameljanczyk A., Długosz P., Strawa M.: Komputerowa implementacja algorytmu wyznaczania wstępnej diagnozy medycznej. In: *VII Konferencja Naukowa Modelowanie Cybernetyczne Systemów Biologicznych, MCSB2010, Kraków*, 2010.
[4]. Andersson J.,Wallace D.: Pareto optimization using the struggle genetic crowding algorithm, *Engineering Optimization*, vol. 34(6), pp. 623–643, 2002.
[5]. Berkhout F., Bouwer L.M., Bayer J., Bouzid M., Cabeza M., Hanger S., Hof A., Hunter P., Meller L., Patt A.: Deep Emissions Reductions and Mainstreaming of Mitigation and Adaptation: Key Findings, 2013.
[6]. Berkhout F., Bouwer L.M., Bayer J., Bouzid M., Cabeza M., Hanger S., Hof A., Hunter P., Meller L., Patt A., et al.: European policy responses to climate change: progress on mainstreaming emissions reduction and adaptation, Regional Environmental Change, vol. 15, pp. 949–959, 2015.
[7]. Christophers S.R.: Aedes aegypti: the yellow fever mosquito, CUP Archive, 1960.
[8]. De Weck O.L.: Multiobjective optimization: History and promise. In: Invited Keynote Paper, GL2-2, The Third China–Japan–Korea Joint Symposium on Optimization of Structural and Mechanical Systems, Kanazawa, Japan, vol. 2, p. 34, 2004.
[9]. Descloux E., Mangeas M., Menkes C.E., Lengaigne M., Leroy A., Tehei T., Guillaumot L., Teurlai M., Gourinat A.C., Benzler J., et al.: Climate-Based Models for Understanding and Forecasting Dengue Epidemics, *PLoS Neglected Tropical Diseases*, vol. 6(2), 2012.
[10]. Ehrgott M.: Vilfredo Pareto and Multi-Objective Optimization, *Documenta Mathematica*, pp. 447–453, 2012.
[11]. Enduri M.K., Jolad S.: Dynamics of dengue disease with human and vector mobility, *Spatial and Spatio-temporal Epidemiology*, vol. 25, pp. 57–66, 2018.
[12]. Erandi K.K.W.H., Perera S.S.N., Mahasinghe A.C.: Dengue Outbreaks Prediction Model for Urban Colombo using Meteorological Data, *International Journal of Dynamical Systems and Differential Equations* (in press).
[13]. He Z., Yen G.G., Zhang J.: Fuzzy-Based Pareto Optimality for Many-Objective Evolutionary Algorithms, *IEEE Transactions on Evolutionary Computation*, vol. 18(2), pp. 269–285, 2013.
[14]. Hii Y.L.: Climate and dengue fever: early warning based on temperature and rainfall. Ph.D. thesis, Umeå University, 2013.
[15]. Gubler D.J.: Epidemic Dengue/Dengue Haemorrhagic Fever: a global public health problem in the 21st century, *Dengue Bulletin*, vol. 12, pp. 1–19, 1997.
[16]. Johansson M.A., Dominici F., Glass G.E.: Local and Global Effects of Climate on Dengue Transmission in Puerto Rico, *PLoS Neglected Tropical Diseases*, vol. 3(2), 2009.
[17]. Karim M.N., Munshi S.U., Anwar N., Alam M.S.: Climatic factors influencing dengue cases in Dhaka city: a model for dengue prediction, *The Indian Journal of Medical Research*, vol. 136(1), p. 32–39, 2012.
[18]. Liu K., Wang T., Yang Z., Huang X., Milinovich G.J., Lu Y., Jing Q., Xia Y., Zhao Z., Yang Y., et al.: Using Baidu Search Index to Predict Dengue Outbreak in China, *Scientific Reports*, vol. 6, p. 38040, 2016.
[19]. Massad E., Ortega N.R.S., de Barros L.C., Struchiner C.J.: Fuzzy Logic in Action: Applications in Epidemiology and Beyond, *Springer Science & Business Media*, 2009.
[20]. McCullagh P., Nelder J.A.: *Generalized Linear Models*, 2nd Edition, Chapman and Hall, London, UK, 1989.
[21]. Morin C.W., Comrie A.C., Ernst K.: Climate and dengue transmission: evidence and implications, *Environmental health perspectives*, vol. 121(11–12), pp. 1264–1272, 2013.
[22]. Ramadona A.L., Lazuardi L., Hii Y.L., Holmner Å., Kusnanto H., Rocklöv J.: Prediction of Dengue Outbreaks Based on Disease Surveillance and Meteorological Data, *PloS one*, vol. 11(3), p. e0152688, 2016.
[23]. Schumpeter J.A.: Vilfredo Pareto (1848–1923), *The Quarterly Journal of Economics*, pp. 147–173, 1949.
[24]. Seidahmed O.M.E., Eltahir E.A.B.: A Sequence of Flushing and Drying of Breeding Habitats of *Aedes Aegypti (L.)* Prior to the Low Dengue Season in Singapore, *PLoS Neglected Tropical Diseases*, vol. 10(7), 2016.
[25]. Wesolowski A., Qureshi T., Boni M.F., Sundsøy P.R., Johansson M.A., Rasheed S.B., Engø-Monsen K., Buckee C.O.: Impact of human mobility on the emergence of dengue epidemics in Pakistan, *Proceedings of the National Academy of Sciences*, vol. 112(38), pp. 11887–11892, 2015.
[26]. Wickramaarachchi W.P.T.M., Perera S.S.N.: Developing a two-dimensional climate risk model for dengue disease transmission in Urban Colombo, *Journal of Basic and Applied Research International*, vol. 20(3), pp. 168–177, 2017.

[27]. Wickramaarachchi W.P.T.M., Perera S.S.N.: The nonlinear dynamics of the dengue mosquito reproduction with respect to climate in urban Colombo: a discrete time density dependent fuzzy model, *International Journal of Mathematical Modelling and Numerical Optimisation*, vol. 8(2), pp. 145–161, 2017.
[28]. Wilder-Smith A.: Dengue vaccine development: status and future, *Bundesgesundheitsblatt Gesundheitsforschung Gesundheitsschutz*, vol. 63(1), pp. 40–44, 2020.
[29]. Wu P.C., Lay J.G., Guo H.R., Lin C.Y., Lung S.C., Su H.J.: Higher temperature and urbanization affect the spatial patterns of dengue fever transmission in subtropical Taiwan, *Science of the Total Environment*, vol. 407(7), pp. 2224–2233, 2009.
[30]. Zadeh L.A.: Toward a theory of fuzzy information granulation and its centrality in human reasoning and fuzzy logic, *Fuzzy Sets and Systems*, vol. 90(2), pp. 111–127, 1997.

**Affiliations**
**Piotr Jakubowski**
Military University of Technology, Faculty of Cybernetics, 00-908 Warsaw, ul. Gen. Kaliskiego 2, Poland, piotr.jakubowski@wat.edu.pl

**Hasitha Erandi**
University of Colombo, Research & Development Center for Mathematical Modeling, Department of Mathematics, Colombo 03, Sri Lanka, hasitha.erandi@yahoo.com

**Anuradha Mahasinghe**
University of Colombo, Research & Development Center for Mathematical Modeling,
Department of Mathematics, Colombo 03, Sri Lanka, anuradhamahasinghe@maths.cmb.ac.lk,
ORCID ID: http://orcid.org/0000-0003-2828-6090

**Sanjeewa Perera**
University of Colombo, Research & Development Center for Mathematical Modeling,
Department of Mathematics, Colombo 03, Sri Lanka, ssnp@maths.cmb.ac.lk

**Andrzej Ameljanczyk**
Military University of Technology, Faculty of Cybernetics, 00-908 Warsaw,
ul. Gen. Kaliskiego 2, Poland, andrzej.ameljanczyk@wat.edu.pl